\documentclass[12pt,a4paper]{article}
\usepackage{amsmath}
\usepackage{amsxtra}
\usepackage{amstext}
\usepackage{amssymb}
\usepackage{latexsym}
\usepackage{dsfont}
\usepackage{graphics}

\newcommand{\nn}{\nonumber}

\newcommand{\p}{\vspace{6pt}\noindent}

\advance\textheight by \topskip

\makeatletter

\@addtoreset{equation}{section}
\def\section{\@startsection {section}{1}{\z@}{-8.5ex plus -1ex minus
 -.2ex}{3.3ex plus .2ex}{\large\bf}}%\centering}}
\def\subsection{\@startsection{subsection}{2}{\z@}{-3.25ex plus
 -1ex minus -.2ex}{1.5ex plus .2ex}{\bf}}
\def\subsubsection{\@startsection{subsubsection}{3}{\z@}{-3.25ex plus%
 -1ex minus -.2ex}{1.5ex plus .2ex}{\sl}}

\begin{document}

\begin{titlepage}
\vspace*{-2cm}
\begin{flushright}
%%YITP-05-29
\end{flushright}

\vspace{0.3cm}

\begin{center}
{\Large {\bf }} \vspace{1cm} {\Large {\bf On the semi-dynamical
reflection equation: solutions and structure matrices}}\\
\vspace{1cm} {\large  J.\ Avan\footnote{\noindent E-mail: {\tt
avan@u-cergy.fr}} and
C.\ Zambon\footnote{\noindent E-mail: {\tt cristina.zambon@u-cergy.fr}} \\
\vspace{1cm}
{\em Laboratoire de Physique Th\'{e}orique et Mod\'{e}lisation \\
Universit\'{e} de Cergy-Pontoise (CNRS UMR 8089), Saint-Martin 2\vspace{.1cm}\\
2 avenue Adolphe Chauvin,
95302 Cergy-Pontoise Cedex, France}}\\

\vspace{1cm} {\bf{ABSTRACT}}
\end{center}
Explicit solutions of the non-constant semi-dynamical reflection
equation are constructed, together with suitable parametrizations of
their structure matrices. Considering the semi-dynamical reflection
equation with rational non-constant Arutyunov-Chekhov-Frolov
structure matrices, and a specific meromorphic ansatz, it is found
that only two sets of the previously found constant solutions are
extendible to the non-constant case. In order to simplify future
constructions of spin-chain Hamiltonians, a parametrization
procedure is applied explicitly to all elements of the
semi-dynamical reflection equation available. Interesting
expressions for `twists' and $R$-matrices entering the
parametrization procedure are found. In particular, some expressions
for the $R$-matrices seem to appear here for the first time. In
addition, a new set of consistent structure matrices for the
semi-dynamical reflection equation is obtained.

\vfill
\end{titlepage}

\section{Introduction}

\p Dynamical extensions of Sklyanin-type quantum reflection 
algebras
\cite{Slyanin88}, \cite{Cherednik84}, \cite{Donin02} have been
introduced and quite extensively studied in the last years
\cite{Arutyunov98}, \cite{Nagy04}. The so called semi-dynamical
reflection algebra, exemplified in \cite{Arutyunov98}, was
formulated generically in \cite{Nagy04} and later used as a basic
algebraic structure to yield formal spin-chain quantum integrable
Hamiltonians \cite{Nagy05}. Generic consistent parametrizations of
its matrices were then proposed in \cite{Avan07}, leading to
remarkably simplified factorized formulas for the generating
monodromy matrices. Together with the classification of scalar
(non-operatorial) solutions $K$ to the specific semi-dynamical
reflection equation (SDRE), started in \cite{Avan05} for the
constant case, this procedure is expected to lead to new, fully
explicit, spin-chain Ruijsenaar-Schneider (RS) type Hamiltonians.
The results presented in this article represent another step in this
direction. Explicit resolutions of the parametrization described in
\cite{Avan07}, together with the subsequent partial classifications
of non-constant scalar $K$ solutions  will be given. Indeed, the
construction of tractable (i.e. locally interacting) spin-chain type
Hamiltonians, contrary to the pure $N$-body system Hamiltonians,
requires a priori the consideration of non-constant, spectral
parameter dependent solutions \cite{Slyanin88}. In this context, the
parametrization proposed in \cite{Avan07} considerably simplifies
the form of the relevant monodromy matrices, and therefore it is a
key ingredient to explicitly build the Hamiltonians.

\p The SDRE is a quadratic constraint equation for generators of a
quantum algebra $\mathcal{G}$ encapsulated into the matrix $K$. Such
a constraint is represented as an equation in
$\mathrm{End}(\mathcal{U}) \otimes \mathrm{End}(\mathcal{U})$ where
$\mathcal{U}$ is a given vector space known as the auxiliary space.
This space can be a finite dimensional space $V$ or a loop space $V
\otimes \mathds{C}[u]$. The general form of the SDRE considered in
this article is:
\begin{equation}\label{SDRE}
A_{12}(u,v;\lambda)\,K_1(u;\lambda)\,B_{12}(v;\lambda)\,K_2(v;\lambda+\gamma
h_1)=K_2(v;\lambda)\,C_{12}(u;\lambda)\,K_1(u;\lambda+\gamma
h_2)\,D_{12}(u-v;\lambda),
\end{equation}
where $A$, $B$, $C$ and $D$ are $\mathds{C}$-number matrices known
as structure matrices. All elements appearing in \eqref{SDRE} depend
on a set of complex variables, collectively denoted $\lambda\equiv
\{\lambda_i,\; i=1,\dots,n\}$ and known as dynamical variables. If
the auxiliary space is a loop space, there is also a spectral
parameter dependence in $A$, $B$, $C$, $D$, $K$ represented by the
complex variables $u$ and $v$. Notice that the matrix $D$ may have a
more general dependence in the spectral parameters, like $A$, but
this leads to subsequent difficulties when deriving the commuting
Hamiltonians \cite{Avan07}. The dynamical variables $\{\lambda_i\}$
are interpreted as coordinates on the dual $\textbf{h}^*$ of a
$n$-dimensional abelian subalgebra $\textbf{h}$ of a simple Lie
algebra $\textbf{g}$.

\p Given a basis $h\equiv \{h^i,\; i=1,\dots,n\}$ of $\textbf{h}^*$
(with $\{h_i,\; i=1,\dots,n\}$ basis in $\textbf{h}$), and setting
$\lambda=\sum_{i=1}^n\, \lambda_i h^i$, it is possible to define
formally
\begin{equation}\label{definitionhplusgamma}
f(\lambda+\gamma h)\equiv e^{\gamma\mathcal{D}}\,f(\lambda)\,
e^{-\gamma\mathcal{D}},\qquad \mathcal{D}=\sum_{i=1}^n\, h_i
\partial_{\lambda_i},
\end{equation}
where $f(\lambda)$ is a differentiable function on $\textbf{h}^*$
and the auxiliary space $\mathcal{U}$ is assumed to be a
diagonalizable irreducible module of $\textbf{h}$. In order to
simplify the notation, it will be often set $f(h)\equiv
f(\lambda+\gamma h)$.

\p The structure matrices are supposed to satisfy the following
zero-weight conditions
\begin{equation}
[h_i \otimes \mathds{1}, B_{12}]=0, \quad [\mathds{1}\otimes h_i,
C_{12}]=0, \quad [h_i \otimes \mathds{1}+\mathds{1}\otimes h_i,
D_{12}]=0, \quad i=1,\dots, n.
\end{equation}
Moreover, the assumed associativity of the semi-dynamical reflection
algebra $\mathcal{G}$ yields, as sufficient consistency conditions,
that the structure matrices obey a set of YB-like equations
\cite{Nagy04}. Such equations have been reformulated in
\cite{Avan07} to take into account certain freedom enjoyed by the
structure matrices due to the form-invariance of the SDRE
\eqref{SDRE} under suitable transformations. In fact, multiplying on
the left hand side equation \eqref{SDRE} by $(g \otimes
\mathds{1})$, where $g$ is an automorphism of the auxiliary space
$\mathcal{U}$, leads to an equivalent - though with a different
definition of structure matrices - formulation of the exchange
relations satisfied by the generators of the algebra $\mathcal{G}$,
which are encapsulated into the unmodified matrix $K$. Taking into
account this property, the YB-like consistency equations can be
written as follows
\begin{eqnarray}\label{ABCDconsitencyequations}
\mbox{\textbf{a}}\; && A_{12}\,A_{13}{^{gg}}\,A_{23}=A_{23}{^{gg}}\,A_{13}\,A_{12}{^{gg}},\nn\\
\mbox{\textbf{b}}\; && A_{12}\,C_{13}{^{g_1}}\,C_{23}=C_{23}{^{g_2}}\,C_{13}\,A_{12}(h_3){^{gg}},\nn\\
\mbox{\textbf{c}}\; && D_{12}\,B_{13}\,B_{23}(h_1){^{g_3}}=B_{23}\,B_{13}(h_2){^{g_3}}\,D_{12},\nn\\
\mbox{\textbf{d}}\; &&
D_{12}(h_3)\,D_{13}\,D_{23}(h_1)=D_{23}\,D_{13}(h_2)\,D_{12},
\end{eqnarray}
where
\begin{equation}\label{shortgnotation}
M_{12}{^{gg}}\equiv g_1 g_2\,M_{12}\,g^{-1}_1 g^{-1}_2, \quad
M_{12}{^{g_1}}\equiv g_1\,M_{12}\,g^{-1}_1,\quad
M_{12}{^{g_2}}\equiv g_2\,M_{12}\,g^{-1}_2,
\end{equation}
and $g$ is the automorphism of the auxiliary space. Note that
(\ref{ABCDconsitencyequations}-\textbf{d}) is the
Gervais-Neveu-Felder (GNF) equation \cite{Gervais84}, and it is
unmodified by this extension. In addition, for consistency
conditions, $B_{12}=C_{21}$.

\p In this article only the sets of matrices $A$, $B$, $C$, $D$ for
which there exists scalar invertible solutions $K$ to \eqref{SDRE}
will be taken into account. At present, only one of these sets is
known, namely the Arutyunov-Chekhov-Frolov (ACF) solution
\cite{Arutyunov98}, which is associated to the RS models based on
the Lie algebra data $gl(n)$ \cite{Ruijsenaars86}. In this case the
structure matrices obey the generalized YB set of equations
\eqref{ABCDconsitencyequations} where the automorphism $g$
represents a shift in the spectral parameter $u$ as
$g=\exp(-\gamma\,\frac{d}{du})$. In the limit of non-spectral
parameter dependence the automorphism reduces to $g=\mathds{1}$. For
this set of matrices $A$, $B$, $C$, $D$, the auxiliary space
$\mathcal{U}$ is a finite dimensional loop space $V \otimes
\mathds{C}[u]$ with $V=\mathds{C}^n$ and $\textbf{h}$ is the Cartan
subalgebra of the Lie algebra $gl(n)$. As a consequence, for the
matrix $K$ in \eqref{SDRE}, the definition
\eqref{definitionhplusgamma} translates into
\begin{equation}
K_2(h_1)\equiv \sum_{j=1}^n\, h_j \otimes K(\lambda_j+\gamma).
\end{equation}
Given a basis $\{e_i,\; i=1,\dots, n\}$ of $V$, $e_{ij}=e_i\otimes
e_j$ (with $e_{jj}\equiv h_j$) represents the usual matrix basis and
the rational ACF structure matrices read \cite{Arutyunov98}
\begin{eqnarray}\label{AwithSP}
A_{12}(u,v;\lambda)&=&A_{12}^{\infty}(\lambda)+\frac{\gamma}{(u-v)}d_{12}+\frac{\gamma}{v}b_{12}-\frac{\gamma}{(u+\gamma)}b_{21},\\
B_{12}(v;\lambda)&=&B_{12}^{\infty}(\lambda)-\frac{\gamma}{(v+\gamma)}b_{12},\qquad C_{12}(u;\lambda)=B_{21}(v;\lambda),\label{BwithSP}\\
D_{12}(u-v;\lambda)&=&D_{12}^{\infty}(\lambda)+\frac{\gamma}{(u-v)}d_{12},\label{DwithSP}
\end{eqnarray}
with
\begin{eqnarray}\label{A}
A_{12}^{\infty}(\lambda)&=&\mathds{1}\otimes \mathds{1}+\sum_{i\neq
j=1}^n\,x_{ij}\,(e_{ii}-e_{ij})\otimes (e_{jj}-e_{ji}),\\
B_{12}^{\infty}(\lambda)&=&\mathds{1}\otimes \mathds{1}+\sum_{i\neq
j=1}^n\,y_{ij}\,e_{jj}\otimes (e_{ii}-e_{ij});\qquad
C_{12}^{\infty}(\lambda)=B^{\infty}_{21}(\lambda),
\label{B}\\
D_{12}^{\infty}(\lambda)&=&\mathds{1}\otimes \mathds{1}+\sum_{i\neq
j=1}^n\,x_{ij}\,(e_{ij} \otimes e_{ji}-e_{ii} \otimes
e_{jj}),\label{D}
\end{eqnarray}
and
\begin{eqnarray}
d_{12}&=&\sum_{i,j=1}^n\,e_{ij}\otimes e_{ji},\qquad
b_{12}=\sum_{i,j=1}^n\,e_{ii}\otimes
e_{ji},\qquad c_{12}\equiv b_{21},\nn\\
x_{ij}&=&\frac{\gamma}{(\lambda_i-\lambda_j)}\equiv\frac{\gamma}{\lambda_{ij}},
\qquad
y_{ij}=\frac{\gamma}{(\lambda_i-\lambda_j-\gamma)}\equiv\frac{\gamma}{(\lambda_{ij}-\gamma)}.
\end{eqnarray}
The matrices (\ref{A}-\ref{D}) will be denoted as rational
`constant' ACF matrices because they do not depend on the spectral
parameter. Similarly, matrices (\ref{AwithSP}-\ref{DwithSP}) will be
called non-constant ACF matrices. They exhibit spectral parameter
dependence and reduce to the previous set (\ref{A}-\ref{D}) in the
limit $u,v,(u-v)\longrightarrow \infty$.

\p The parametrization procedure proposed in \cite{Avan07} makes use
of quantum group-like objects such as $R$-matrices and Drinfeld's
twists for building the $A$, $B$, $C$, $D$ and $K$ matrices. It
allows to simplify significantly the expressions for the monodromy
matrices found previously \cite{Nagy05}, and therefore to facilitate
the explicit construction of integrable spin-chain Hamiltonians.

\p The purpose of the present article is to provide explicit
realizations of the parametrization procedure proposed in
\cite{Avan07} making use of the specific matrices $D$ available in
\eqref{DwithSP} and \eqref{D}, and to classify, at least partially,
non-constant solutions $K$ of the SDRE. The first part of the
article is focused on the search of solutions of equation
\eqref{SDRE} using the structure matrices
(\ref{AwithSP}-\ref{DwithSP}). The final aim is to extend the
results obtained in \cite{Avan05} for the rational constant
structure matrices (\ref{A}-\ref{D}) to the case with spectral
parameter dependence. The second part of the paper will be devoted
to the parametrization procedure for which three distinct situations
will be considered. Making use of the ACF set of structure matrices,
the parametrization procedure will be applied first to the simpler
case $g=\mathds{1}$ and later to the more complicated situation in
which $g=\exp(-\gamma\,\frac{d}{du})$. Finally, adopting an
alternative parametrization for matrices $D$ in \eqref{D} and
\eqref{DwithSP}, which is provided in \cite{Antonov97}, it will be
shown how the parametrization procedure leads to an alternative set
of solutions of equations \eqref{ABCDconsitencyequations}, namely to
new structure matrices $A$, $B$, $C$. Full analysis of the SDRE
built by these matrices will be left to further studies, even if
some information concerning the solutions of this equation can be
deduced making use of the parametrization procedure. It should be
emphasized that the existence of this new set of structure matrices
relies on the availability of distinct de-twisting procedures for a
single $D$-matrix. In fact, from \cite{Avan07} a set of consistent
structure matrices $A$, $B$, $C$, $D$ is provided as soon as a
$D$-matrix is chosen and a particular cocycle twist-like formulation
of $D$ is specified. Further remarks concerning this point will be
added later.

\section{Solutions of the non constant rational SDRE}
\label{solutionsofSDRE}

\p In this section, making use of the rational ACF set of solutions
for the equations \eqref{ABCDconsitencyequations}, the problem to
discover and classify the matrices $K$ solving the SDRE will be
addressed. In \cite{Avan05}, as it was pointed out before, this
problem has been already tackled in the case of no spectral
parameter dependence. Four sets of solutions were identified, namely

\begin{eqnarray}\label{KsolutionswoSP}
\mbox{\textbf{Ia}}\;&
&K^{\infty}_{ij}(\lambda)=\frac{f+\lambda_i}{f+\lambda_j}\prod_{a\neq
j}^n\,\frac{\gamma}{\lambda_{ja}}, \nn\\
\mbox{\textbf{IIa}}\;&
&K^{\infty}_{ij}(\lambda)=(f-\Lambda_{ij})\prod_{a\neq
j}^n\,\frac{\gamma}{\lambda_{ja}}, \qquad \Lambda_{ij}=\sum_{a=1}^n\,\lambda_a-(\lambda_i+\lambda_j)\nn\\
\mbox{\textbf{Ib}}\;&
&K^{\infty}_{ij}(\lambda)=\frac{f+\Lambda_j}{f+\Lambda_i}\prod_{a\neq
j}^n\,\frac{\gamma}{\lambda_{ja}}, \qquad \Lambda_i=\sum_{a=1}^n\,\lambda_a-(\lambda_i)\nn\\
\mbox{\textbf{IIb}}\;&
&K^{\infty}_{ij}(\lambda)=\frac{f}{f-\lambda_{ij}}\prod_{a\neq
j}^n\,\left(1+\frac{f}{\lambda_{ja}}\right),
\end{eqnarray}
where $f$ is a function $\gamma$-periodic on each dynamical
variable. All other solutions can be obtained from these sets by
multiplying the $K$-matrix on the right by a diagonal matrix
$N_{ii}(\lambda)$ satisfying the following flatness condition
\begin{equation}
N_{ii}(\lambda)\,N_{jj}(\lambda+\gamma
h_i)=N_{jj}(\lambda)\,N_{ii}(\lambda+\gamma h_j).
\end{equation}
In particular, the matrix $N_{jj}=0, \, N_{ii}=1, \, i\neq j$ allows
to obtain a matrix $K$ with the entries $K_{kj}=0,\, k=1,\dots,n$.
In this case, the periodicity condition on the dynamical variable
$\lambda_j$ for the function $f$ is omitted. It can be noticed that
for $n=2$ solutions (\ref{KsolutionswoSP}-\textbf{Ia}) and
(\ref{KsolutionswoSP}-\textbf{IIa}) collapse to
(\ref{KsolutionswoSP}-\textbf{Ib}) and
(\ref{KsolutionswoSP}-\textbf{IIb}), respectively. Moreover, the
only invertible matrices are represented by solutions
(\ref{KsolutionswoSP}-\textbf{IIb}), since matrices
(\ref{KsolutionswoSP}-\textbf{Ia}) and
(\ref{KsolutionswoSP}-\textbf{Ib}) have rank 1 while matrices
(\ref{KsolutionswoSP}-\textbf{IIa}) have rank 2. Finally, note that
in the limit $f\longrightarrow 0$ solutions
(\ref{KsolutionswoSP}-\textbf{IIb}) reduce to the trivial solution
$\mathds{1}$.

\p The classification of the matrices $K$ solving the SDRE will now
be extended to the case with spectral parameter dependence. In order
to simplify the notation, the explicit dependence on the dynamical
variables will be omitted in what follows. Because of the particular
form of the structure matrices, and without making any assumption on
the unknown matrices $K$, the SDRE \eqref{SDRE} can be rewritten in
a more appealing form as
\begin{eqnarray}\label{SDREsimplified}
&&\left(A_{12}(u,v)-\frac{\gamma}{(u-v)}d_{12}
\right)\,K_1(u)\,B_{12}(v)\,K_2(v; h_1)\nn\\
&=&\left(\frac{}{}K_2(v)\,C_{12}(u)\,K_1(u;h_2)
-K_2(u)\,C_{12}(v)\,K_1(v; h_2)\right)\frac{\gamma}{(u-v)}d_{12}\nn\\
&&+K_2(v)\,C_{12}(u)\,K_1(u; h_2)\,D^{\infty}_{12}.
\end{eqnarray}
The advantage of this formulation is to gather together in a more
compact way the terms proportional to the factor $1/(u-v)$. Then,
the following ansatz for the matrix $K$ will be used
\begin{equation}\label{Kansatz}
K(u,\lambda)=\sum_{l=0}^{N}\left(\frac{\gamma }{u}\right)^l
k^{(l)}(\lambda),\qquad k^{(0)}(\lambda)\equiv K^{\infty}(\lambda).
\end{equation}
This expansion in powers of $1/u$ represents a natural extension, as
rational function, of the solutions \eqref{KsolutionswoSP} to which
it reduces in the limit $u\longrightarrow\infty$. Note that in
\eqref{Kansatz} the location of the poles at $u=0$ is just a matter
of choice. Multiplying such an ansatz by
\begin{equation}\label{mulfactor}
\prod_{l=1}^{N}\,\left(\frac{u}{u-u_0}\right)^l
\end{equation}
allows to obtain an equivalent ansatz with the poles shifted at
$u=u_0$. Reciprocally, given any matrix $K$ with a finite set of
poles, it can always be brought back to the form \eqref{Kansatz} by
a suitable multiplicative factor like \eqref{mulfactor}.

\p Once \eqref{Kansatz} is plugged into the expression
\eqref{SDREsimplified}, it is noticed that all terms coming from the
second line, and proportional to $1/(u-v)$ can be combined together
in such a way to eliminate completely this factor. For instance:
\begin{equation}
\frac{1}{(u-v)}\left[\frac{1}{
v(u+\gamma)}-\frac{1}{u(v+\gamma)}\right]=\frac{\gamma}{u
v(u+\gamma)(v+\gamma)}.
\end{equation}
Then, using the property
\begin{equation}
c_{12}\,k^{(l)}_1\,b_{12}\,k^{(l)}_2(h_1)=k^{(l)}_2\,c_{12}\,k^{(l)}_1(h_2)\,d_{12},\qquad
l=0,1,\dots, N
\end{equation}
some simplifications can be performed amongst terms coming from the
first and the second lines of \eqref{SDREsimplified}. The remaining
terms must be treated with care. First of all, the powers of the
spectral parameters appearing in each term must be reduced as much
as possible by decomposition in prime elements. Finally, making use
of the property
\begin{equation}
(A^{\infty}-b)_{12}\,k^{(l)}_1b_{12}\,k^{(l)}_2(h_1)=0,\qquad
l=0,1,\dots, N
\end{equation}
and the reduction explained above, further simplifications are
possible. The expression obtained splits into several relations,
each of them gathering algebraically independent terms. They
represents constraints for the elements of the matrices $k^{(l)}$,
which must be analyzed carefully by projecting them onto the matrix
elements $(e_{ij} \otimes e_{kl})\;\; i,j,k,l=1,\dots,n$.

\p Starting with the simplest ansatz for the matrices $K$, namely
the expression \eqref{Kansatz} with $l=1$, the expression
\eqref{SDREsimplified} translates into $8$ relations, namely

\begin{eqnarray}
&&A^{\infty}_{12}\,k^{(0)}_1B^{\infty}_{12}\,k^{(0)}_2(h_1)
=k^{(0)}_2\,C^{\infty}_{12}\,k^{(0)}_1(h_2)\,D^{\infty}_{12},\label{inftyconstraint}\\
&&\nn\\
&&
b_{12}\,k^{(0)}_1\,(B^{\infty}-b)_{12}\,k^{(1)}_2(h_1)=0,\label{firstconstraint}\\
&&\nn\\
&&b_{12}\,k^{(1)}_1\,(B^{\infty}-b)_{12}\,k^{(1)}_2(h_1)=0,\nn\\
&&(A^{\infty}-c)_{12}\,k^{(1)}_1\,B^{\infty}_{12}\,k^{(1)}_2(h_1)
-k^{(1)}_2\,(C^{\infty}-c)_{12}\,k^{(1)}_1(h_2)\,D^{\infty}_{12}\nn\\
&&+k^{(0)}_2\,(C^{\infty}-c)_{12}\,k^{(1)}_1(h_2)\,d_{12}+(b-d)_{12}
\,k^{(1)}_1(B^{\infty}-b)_{12}\,k^{(0)}_2(h_1)=0,\nn\\
&&(A^{\infty}-c)_{12}\,k^{(1)}_1\,B^{\infty}_{12}\,k^{(0)}_2(h_1)
-k^{(0)}_2\,(C^{\infty}-c)_{12}\,k^{(1)}_1(h_2)\,D^{\infty}_{12}=0,\label{zeroconstraint}\\
&&\nn\\
\mbox{\textbf{a}}\;&
&c_{12}\,(k^{(1)}-k^{(0)})_1\,B^{\infty}_{12}\,k^{(0)}_2(h_1)
-k^{(0)}_2\,c_{12}\,(k^{(1)}-k^{(0)})_1(h_2)\,D^{\infty}_{12}=0,\nn\\
\mbox{\textbf{b}}\;& &c_{12}\,(k^{(1)}-k^{(0)})_1\,(B^{\infty}_{12}-b)_{12}\,k^{(1)}_2(h_1)
-k^{(1)}_2\,c_{12}\,(k^{(1)}-k^{(0)})_1(h_2)\,(D^{\infty}_{12}-d)_{12}=0,\nn\\
\mbox{\textbf{c}}\;&
&A^{\infty}_{12}\,k^{(0)}_1\,B^{\infty}_{12}\,k^{(1)}_2(h_1)-k^{(1)}_2\,C^{\infty}_{12}\,k^{(0)}_1(
h_2)\,D^{\infty}_{12}+
b_{12}\,k^{(0)}_1\,(B^{\infty}-b)_{12}\,k^{(0)}_2(h_1)=0.\label{Kconstraints}
\end{eqnarray}
Equation \eqref{inftyconstraint}, which involves only the matrix
$k^{(0)}$, has already been investigated and its solutions have been
listed in \eqref{KsolutionswoSP} ($k^{(0)}\equiv K^{\infty}$). All
other  8 relations incorporate both $k^{(0)}$ and $k^{(1)}$
matrices. When they are analyzed, one discovers that
\eqref{firstconstraint} represents a strong constraint for the
matrix $k^{(1)}$.

\p Consider first the case when $k^{(0)}$ has no zero entry, then
expression \eqref{firstconstraint} states the following
\begin{equation}
k^{(1)}_{ij}=k^{(1)}_{kj}, \qquad i\neq k=1,\dots,n.
\end{equation}
As a consequence the relations \eqref{zeroconstraint} become
identities and only the 3 relations \eqref{Kconstraints} remain to
be investigated. For instance, (\ref{Kconstraints}-\textbf{a})
allows to establish whether the sets of solutions listed in
\eqref{KsolutionswoSP} can be extended or not, and to specify the
form of the extensions. Provided $k^{(1)}$ has no zero entry, it
turns out that (\ref{KsolutionswoSP}-\textbf{Ia}) and
(\ref{KsolutionswoSP}-\textbf{IIa}) are not extendable to a first
order solution $k^{(1)}$. On the contrary, the solutions
(\ref{KsolutionswoSP}-\textbf{Ib}) and
(\ref{KsolutionswoSP}-\textbf{IIb}) can be uniquely extended as
follows
\begin{eqnarray}\label{Ksolutions}
\mbox{\textbf{Ib}}\;&
&K_{ij}(u;\lambda)=(f+\Lambda_j)\left(\frac{1}{f+\Lambda_i}-\frac{1}{u}\right)\prod_{a\neq
j}\;\left(\frac{\gamma}{\lambda_{ja}}\right),
\qquad \Lambda_i=\sum_{a=1}^n\,\lambda_a-(\lambda_i),\nn\\
\mbox{\textbf{IIb}}\;&
&K_{ij}(u;\lambda)=f\left(\frac{1}{f-\lambda_{ij}}-\frac{1}{u}\right)\prod_{a\neq
j}\;\left(\frac{f}{\lambda_{ja}}+1\right).
\end{eqnarray}
The constraints provided by expressions
(\ref{Kconstraints}-\textbf{b}) and (\ref{Kconstraints}-\textbf{c})
are automatically satisfied by the sets of solutions
\eqref{Ksolutions}. Notice that solution
(\ref{Ksolutions}-\textbf{IIb}) with constant $f$ coincides with the
solution found by ACF in \cite{Arutyunov98}. Moreover, it can be
shown that allowing a column in the matrix $k^{(1)}$ to be zero,
constraints \eqref{Kconstraints} force all elements of the
corresponding column in $k^{(0)}$ to be identical. However, this can
only happen provided these elements are equal to zero, which is not
allowed by the starting hypothesis. Then, the case $k^{(0)}$ with no
zero entry is completely covered.

\p Consider now the case when $k^{(0)}$ has a zero entry, and
therefore the whole column is zero, as established in \cite{Avan05}.
Solutions $K^{\infty}=k^{(0)}$ with one or more zero-columns can be
obtained by performing suitable simple transformations on the full
solutions \eqref{KsolutionswoSP}, and this then builds the whole set
of constant solutions. Here similar results can be established.
Assuming $k^{(0)}$ has a column set to zero, the corresponding
column in $k^{(1)}$ is forced to be zero as well. Therefore, it is
possible to conclude that solutions $K$ with columns set to zero are
possible, and they are obtained by setting to zero one or several
columns of solutions \eqref{Ksolutions}, since the specific form of
non-zero columns does not depend on the existence of other
zero-columns. Finally, it can be noticed that the limit
$f\longrightarrow 0$ in (\ref{Ksolutions}-\textbf{IIb}) provides
again the solution $\mathds{1}$ where the spectral parameter does
not appear. However, such a solution can be extended in a way to
include a spectral parameter dependence as follows
\begin{equation}\label{diagsolutions}
K(u;\lambda)=\mathds{1}\left(1+\frac{f'}{u}\right),
\end{equation}
where $f'$ is any function $\gamma$-periodic on each dynamical
variable.

\p Attempts to find alternative solutions by truncating the
expansion \eqref{Kansatz} to orders higher than $l=1$ proved to be
unsuccessful. First of all, it is possible to show that to the order
$l=2$ the ansatz \eqref{Kansatz} with $k^{(0)}$ and $k^{(1)}$ given
by \eqref{Ksolutions} is not a solution of the set of equations
coming from expression \eqref{SDREsimplified}, unless $k^{(2)}$=0.
Therefore, the only way out is to reconsider the situation with both
matrices $k^{(1)}$ and $k^{(2)}$ unknown. In fact, the relation
\eqref{inftyconstraint} is the only constraint which emerged
unaltered by using the ansatz \eqref{Kansatz} for a generic order
$l$. All other relations stemming from \eqref{SDREsimplified} depend
on the order of the ansatz \eqref{Kansatz} chosen. For the order
$l=2$ this possibility has been analyzed in detail for the two
solutions (\ref{KsolutionswoSP}-\textbf{Ia}) and
(\ref{KsolutionswoSP}-\textbf{IIb}). It is found that no non trivial
extensions matching the ansatz \eqref{Kansatz} with $l=2$ are
allowed. This lack of success for the ansatz \eqref{Kansatz} with
$l=2$ suggests similar conclusions also hold for an ansatz with a
higher value of $l$.

\section{Parametrization procedure for the elements of the SDRE}
\label{parametrization}

\p In this section the principal formulae of the parametrization
proposed in \cite{Avan07} and obtained by solving equations
\eqref{ABCDconsitencyequations} will be summarized. In what follows,
the dependence of spectral parameters and dynamical variables is
implicit. When a quantity is non-dynamical it will be clearly
stated. It is also assumed, as usual, that matrices $A$, $B$, $C$,
$D$ are invertible.

\p First, the equation (\ref{ABCDconsitencyequations}-\textbf{c})
together with the fact that $B$ is a space-$1$ zero weight matrix,
allows to establish the existence of an invertible $(n\times n)$
matrix $b$ such that
\begin{equation}\label{paramB}
B_{12}\equiv C_{21}=b_2^{-1}\,b_2{^{g}}(h_1).
\end{equation}
This parametrization for the matrix $B$, by means of equation
(\ref{ABCDconsitencyequations}-\textbf{b}), allows to prove the
existence of a quasi non-dynamical $R$-matrix such that
\begin{equation}\label{paramA}
A_{12}=b_1^{-1}\,(b_2{^{g}}){^{-1}}\,R_{12}\,(b_1){^{g}}\,b_2.
\end{equation}
The $R$-matrix appearing in \eqref{paramA} is said to be quasi
non-dynamical since the occurrence of the automorphism $g$ in
(\ref{ABCDconsitencyequations}-\textbf{b}) leads to the following
constraint
\begin{equation}\label{nondynproperty}
R_{12}=R_{12}{^{gg}}(h_3).
\end{equation}
As a consequence
\begin{equation}\label{Rzero}
R_{12}=(e^{-\frac{\sigma}{\gamma}(\log g_1+\log
g_2)})\,R^0_{12}\,(e^{\frac{\sigma}{\gamma}(\log g_1+\log
g_2)}),\qquad \sigma=\sum_{k=1}^n\,\lambda_k
\end{equation}
where the matrix $R^0$, by means of equation
(\ref{ABCDconsitencyequations}-\textbf{a}), is proved to be a
non-dynamical solution of the following $g$-deformed YBE
\begin{equation}\label{shiftedYBE}
R^0_{12}R^0_{13}{^{gg}}R^0_{23}=R^0_{23}{^{gg}}R^0_{13}R^0_{12}{^{gg}}.
\end{equation}
Finally, the matrix $D$ is assumed to be decomposable as
\begin{equation}\label{parametrizationforD}
D_{12}=q_1^{-1}(h_2)\,q_2 ^{-1}\,\tilde{R}_{12}\,q_1\, q_2(h_1),
\end{equation}
where $\tilde{R}$ is also of the form \eqref{Rzero}. It should be
pointed out that all constant $D$-matrices of weak Hecke type
\cite{Faddeev90} associated to a Lie algebra $\textbf{g}=gl(n)$
($n\geq 2$) admit such a decomposition \cite{Buffenoir06} with
$g=\mathds{1}$. This result was recently extended to the affine
(trigonometric) dynamical $R$-matrices \cite{Buffenoir07}. Such
decompositions were already known in a number of cases (see for
example \cite{Antonov97} \cite{Arnaudon98}). Indeed, they
characterize the matrix $D$ as representation of a particular
(cocycle) Drinfeld's twist \cite{Drinfeld90} acting on a universal
$R$-matrix, to yield a quasi-Hopf algebra structure. Note that one
can show immediately the following \cite{Babelon}:

\p \textbf{Proposition 3.1}

\p If $R$ obeys the YBE \eqref{shiftedYBE}, and $D$, which is
constructed from $R$ as \eqref{Rzero} and
\eqref{parametrizationforD} for some $q$, is a zero-weight matrix,
namely $[h_i \otimes \mathds{1}+\mathds{1}\otimes h_i, D_{12}]=0$
for $i=1,\dots, n$, then $D$ obeys the GNF equation. In other words,
the zero-weight condition is sufficient in the cocycle formulation
\eqref{parametrizationforD}.

\p Using parametrization \eqref{paramB}-\eqref{parametrizationforD}
for structure matrices, one finds a consistent (sufficient)
parametrization for the scalar solutions $K$ of the corresponding
SDRE. In particular, one finds as a solution
\begin{equation}\label{paramKbq}
K=(b_{\phantom{}}{^{g}}){^{-1}}\,Q\,q,
\end{equation}
where $Q$ solves
\begin{equation}\label{RQequation}
R_{12}\,Q_1\,q_1\,Q_2(h_1)\,q_1^{-1}
=Q_2\,q_2\,Q_1(h_2)\,q_2^{-1}\,\tilde{R}_{12}.
\end{equation}
Since in all situations analyzed in the present article
$K=\mathds{1}$ is a solution of the SDRE, it is easily shown that
one can choose $R=\tilde{R}$ and $q=b_{\phantom{}}{^g}$ in
\eqref{parametrizationforD}. Therefore, consider the following
parametrization for the matrix $D$
\begin{equation}\label{paramD}
D_{12}=(b_1{^{g}}){^{-1}}(h_2)\,(b_2
{^{g}}){^{-1}}\,R_{12}\,b_1{^{g}}\, b_2{^{g}}(h_1).
\end{equation}
If $Q$ is searched for as quasi non-dynamical, namely
\begin{equation}\label{paramKbb}
Q\equiv e^{-\frac{\sigma}{\gamma}\log
g}Q^0\,e^{\frac{\sigma}{\gamma}\log g},
\end{equation}
with $Q^0$ non-dynamical, then \eqref{RQequation} simplifies to the
following modified YB-like equation \footnote{Note that a more
general situation is represented by replacing $g$ with $\tilde{g}$
in \eqref{paramKbb} and \eqref{RQzeroequation} provided
$[g,\tilde{g}]=0.$}
\begin{equation}\label{RQzeroequation}
R^0_{12}\,Q^0_1\,g_2^{-1}\,Q^0_2\,g_2
=Q^0_2\,g_1^{-1}\,Q^0_1\,g_1\,R^0_{12}.
\end{equation}
This description emphasizes that a systematic scheme to build $A$,
$B$, $C$, $D$ structure matrices arises. Indeed, starting from a
given $D$-matrix, with a specific decomposition
\eqref{parametrizationforD} yielding an $\tilde{R}$-matrix
($\tilde{R}=R$) with a quasi non-dynamical property
\eqref{nondynproperty} such that the associated non-dynamical matrix
$R^0$ obeys a $g$-deformed YB equation \eqref{shiftedYBE}, one has
all the ingredients - namely $R$, $b$ and $g$ - to consistently
build the remaining $A$, $B$, $C$ matrices.

\p At this stage it is worth recalling that non-scalar, operatorial
solutions to the SDRE can also be obtained from
\eqref{RQzeroequation} and \eqref{paramKbq}. In fact, it is
remarkably simple to prove the transfer matrix formula for the SDRE,
such as obtained in \cite{Nagy04}, when it is expressed under a
factorized form (see \cite{Avan07}). Indeed, it is possible to show
the general 'dynamization of trace` as follows:

\p \textbf{Proposition 3.2}

\p Suppose $Q^0$ is a non-dynamical representation, for instance a
monodromy matrix, of \eqref{RQzeroequation} for a given $R$-matrix
$R^0$, and suppose $b$ is a dynamical matrix in End($\mathcal{U}$)
such that
$D_{12}=b_1^{-1}(h_2)\,b_2 ^{-1}\,R^0_{12}\,b_1\, b_2(h_1)$ is
zero-weight (we recall that the auxiliary space is $\mathcal{U}=V \otimes
\mathds{C}[u]$, and we take $g=\mathds{1}$ for simplicity). It is then
possible to construct a dynamical transfer matrix
$\tau^0=(b^{-1}\,Q^0\,b\,e^{\partial_{\lambda}})$ such that
$$[\mathrm{Tr}_V \,\tau^0(u),\mathrm{Tr}_V\, \tau^0(v)]=0.$$
\emph{Proof.} It is simple, provided the following technical tricks
are used:

\p (a) for any three operators $\mathcal{M}_1$, $\mathcal{N}_2$,
$\mathcal{O}_2$ not containing $e^{\partial_{\lambda}}$
$$\mathrm{Tr}_{12}(\mathcal{M}_1\,
e^{\partial_1}\,\mathcal{N}_2\,\mathcal{O}_2\,e^{\partial_2})=
\mathrm{Tr}_{12}(\mathcal{N}_2(h_1)\,\mathcal{M}_1\,e^{\partial_1}\,\mathcal{O}_2\,e^{\partial_2}),$$

\p (b) for any zero-weight $\mathds{C}$-number matrix $D_{12}$ and
any operator $\mathcal{O}_{12}$ not containing
$e^{\partial_{\lambda}}$
$$\mathrm{Tr}_{12}(D_{12}\,\mathcal{O}_{12}\,D^{-1}_{12}\,e^{(\partial_{1}+\partial_{2})})
=\mathrm{Tr}_{12}(\mathcal{O}_{12}\,e^{(\partial_{1}+\partial_{2})}).\footnote{Two
cases of this `dynamical cyclicity' appear in \cite{Doikou05}}$$
Therefore, dynamical trace formula $\tau^0$ seems to be the
one-space counterpart of the dynamical cocycle formula
\eqref{parametrizationforD}.

\subsection{Parametrization for the ACF rational constant matrices}

\p The parametrization procedure will first be applied as an
exercise to the structure matrices (\ref{A}-\ref{D}). In this case
the auxiliary space $\mathcal{U}$ reduces to a finite dimensional
vector space $V$ and $g=\mathds{1}$.

\p Using expression \eqref{B} for $B$, a solution $b$ of
\eqref{paramB} is
\begin{equation}\label{binfty}
b^{\infty}_{ij}=\prod_{1=a\neq
j}^n\,\frac{\lambda_j^{(i-1)}}{\lambda_{ja}}\quad i,j=1\dots n,
\end{equation}
where the notation $b^{\infty}$ emphasizes the independence from a
spectral parameter. Note that such a solution is not unique. In
fact, alternative solutions can be obtained by multiplying each row
$i$ of the matrix \eqref{binfty} by a function $f_i$
$\gamma$-periodic on each dynamical variable. Furthermore, it should
be kept in mind that a matrix obtained by interchanging each pair of
rows in \eqref{binfty} is still a solution of \eqref{paramB}.

\p Similarly, a non-dynamical $R$-matrix solving \eqref{paramA}
reads
\begin{equation}\label{Rinfty}
R^{\infty}=\sum_{i,j=1}^n(e_{ii}\otimes
e_{jj})+\gamma\sum_{i=1}^n\sum_{k=1}^{i-1}\sum_{j=1}^{i-k}(e_{ii-k}\otimes
e_{jj+k-1}-e_{jj+k-1}\otimes e_{ii-k}),
\end{equation}
which is a Cremmer-Gervais $R$-matrix type \cite{Cremmer90}. Note
that in this rational case there is a difference of one in the root
heights, compared to the trigonometric case \cite{Buffenoir06}.
Clearly, in this case $R=R^0\equiv R^{\infty}$ and the matrix
$R^{\infty}$ solves directly the ordinary YBE. The consequent
explicit parametrization of the matrix $D^{\infty}$ \eqref{paramD}
in terms of the non-dynamical $R^{\infty}$-matrix \eqref{Rinfty} and
the `twist' matrix $b^{\infty}$ \eqref{binfty} provides a concrete
example of the theorem mentioned previously and proved in
\cite{Buffenoir06}.

\p Finally, concerning the constant solutions $K$ of the SDRE, it
can be verified that their corresponding $Q=Q^{\infty}$ matrices
satisfy equation \eqref{paramKbb}. Particularly, for solutions
(\ref{KsolutionswoSP}-\textbf{IIa}) and
(\ref{KsolutionswoSP}-\textbf{IIb}), these matrices are non
dynamical and their expressions are particularly simple. For them,
the $Q^{\infty}=Q^0$ matrices are
\begin{eqnarray}\label{Qinfty}
&&\mbox{\textbf{IIa}}\qquad Q^{\infty}=(f\,e_{n n}+e_{n-1 1}+ e_{n
2})\,\gamma^{(n-1)},\nn\\
&&\nn\\
&&\mbox{\textbf{IIb}}\qquad
Q^{\infty}=\mathds{1}+\sum_{i>j=1}^{n}\,\left(
                                          \begin{array}{c}
                                            i-1 \\
                                            i-j \\
                                          \end{array}
                                        \right)
f^{(i-j)}\,e_{ij}.
\end{eqnarray}
As expected, (\ref{Qinfty}-\textbf{IIa}) is a set of rank $2$
matrices, while (\ref{Qinfty}-\textbf{IIb}) is an invertible
triangular set of matrices. By contrast, solutions
(\ref{KsolutionswoSP}-\textbf{Ia}) and
(\ref{KsolutionswoSP}-\textbf{Ib}) cannot be de-dynamized by
\eqref{paramKbq}.

\subsection{Parametrization for the ACF rational non-constant matrices}
\label{paramewithSP}

\p In this situation the structure matrices
(\ref{AwithSP}-\ref{DwithSP}) are solutions of equation
\eqref{ABCDconsitencyequations} with
$g=\exp(-\gamma\,\frac{d}{du})$. The parametrization presented
formally in section (\textbf{\ref{parametrization}}) is still
available, even if it appears to be a little more cumbersome. The
`twist' $b$ matrix is chosen as
\begin{equation}\label{dynamicaltwistwithsp}
b_{ij}=\prod_{a\neq j}\;\frac{\lambda_j^{(i-1)}}{\lambda_{ja}},\quad
b_{nj}=\prod_{a\neq
j}\;\frac{\lambda_j^{(n-1)}(\sigma-\lambda_j+u+\gamma)}{\lambda_{ja}(\sigma+u+f_0)}\quad
j=1\dots n,\quad i=1\dots n-1,
\end{equation}
where $f_0$ is a function of $(\sigma+u)$. For simplicity, from now
on, it will be taken to be a constant. Note that the spectral
parameter dependence is limited to one single row of the matrix.
Even this time this solution is not unique. Instead, the $R$-matrix
reads
\begin{eqnarray}\label{dynamicalRmatrixsln}
R&=&\left(1+\frac{\gamma}{u-v}\right)\sum_{i=1}^{n-1}e_{ii}\otimes
e_{ii}+\left(1+\frac{\gamma}{u-v}\right)\frac{(\sigma+v+f_0)(\sigma+u-\gamma+f_0)}{(\sigma+u+f_0)(\sigma+v-\gamma+f_0)}\;e_{nn}\otimes
e_{nn}\nn\\
&&+\sum_{ij=1}^{n-1}e_{ii}\otimes
e_{jj}+\frac{(\sigma+v+f_0)}{(\sigma+v-\gamma+f_0)}\sum_{i=1}^{n-1}e_{ii}\otimes
e_{nn}+\frac{(\sigma+u-\gamma+f_0)}{(\sigma+u+f_0)}\sum_{i=1}^{n-1}e_{nn}\otimes
e_{ii}\nn\\
&&+\gamma\frac{(\sigma+u+\gamma)}{(\sigma+u+f_0)}\sum_{k=1}^{n-1}\sum_{j=1}^{n-k}e_{ii-k}\otimes
e_{jj+k-1}
-\gamma\frac{(\sigma+v)}{(\sigma+v-\gamma+f_0)}\sum_{k=1}^{n-1}
\sum_{j=1}^{n-k}e_{jj+k-1}\otimes
e_{ii-k}\nn\\
&&+\gamma\sum_{i=1}^{n-1}\sum_{k=1}^{i-1}\sum_{j=1}^{i-k}(e_{ii-k}\otimes
e_{jj+k-1}-e_{jj+k-1}\otimes
e_{ii-k})+\left(\frac{\gamma}{u-v}\right)\sum_{i\neq
j=1}^{n-1}e_{ij}\otimes
e_{ji}\nn\\
&&+\left(\frac{\gamma}{u-v}\right)\frac{(\sigma+u-\gamma+f_0)}{(\sigma+v-\gamma+f_0)}\sum_{i=1}^{n-1}e_{in}\otimes
e_{ni}+\left(\frac{\gamma}{u-v}\right)\frac{(\sigma+v+f_0)}{(\sigma+u+f_0)}\sum_{i=1}^{n-1}e_{ni}\otimes
e_{in}\nn\\
&&-\frac{\gamma}{(\sigma+u+f_0)}\sum_{i=1}^{n-1}\sum_{k=1}^{n-1-i}e_{nn-i}\otimes
e_{kk+i}+\frac{\gamma}{(\sigma+v-\gamma+f_0)}\sum_{i=1}^{n-1}\sum_{k=1}^{n-1-i}e_{kk+i}\otimes
e_{nn-i}.\nn\\
\end{eqnarray}
Note that the results obtained in the previous section are
reproduced when the spectral parameters goes to infinity and
consequently expressions \eqref{dynamicaltwistwithsp} and
\eqref{dynamicalRmatrixsln} reduce to \eqref{binfty} and
\eqref{Rinfty}, respectively. As an example, for the specific case
$n=2$ expression \eqref{dynamicalRmatrixsln} becomes
\begin{eqnarray}\label{dynamicalRmatrixsl2}
R&=&\left(1+\frac{\gamma}{u-v}\right)\left(e_{11}\otimes
e_{11}+\frac{(\sigma+v+f_0)(\sigma+u-\gamma+f_0)}{(\sigma+v-\gamma+f_0)(\sigma+u+f_0)}
\;e_{22}\otimes
e_{22}\right)\nn\\
&&+\frac{(\sigma+v+f_0)}{(\sigma+v-\gamma+f_0)}\;e_{11}\otimes
e_{22}+\frac{(\sigma+u-\gamma+f_0)}{(\sigma+u+f_0)}\;e_{22}\otimes
e_{11}\nn\\
&&-\gamma\frac{(\sigma+v)}{(\sigma+v-\gamma+f_0)}\;e_{11}\otimes
e_{21}+\gamma\frac{(\sigma+u+\gamma)}{(\sigma+u+f_0)}\;e_{21}\otimes
e_{11}\nn\\
&&+\left(\frac{\gamma}{u-v}\right)\left(\frac{\sigma+u-\gamma+f_0}{\sigma+v-\gamma+f_0}\;e_{12}\otimes
e_{21}+\frac{\sigma+v+f_0}{\sigma+u+f_0}\;e_{21}\otimes
e_{12}\right).
\end{eqnarray}
Notice that the choice $f_0=\gamma$ enables to simplify a little
this expression for the $R$-matrix.

\p As expected, and unlike the previous case, the $R$-matrix is
still dynamical. However, in agreement with \eqref{Rzero}, the
dynamical dependence can indeed be eliminated and, setting
$f_0=\gamma$, the matrix $R^0$ reads
\begin{eqnarray}\label{R0matrixsln}
R^0&=&\left(1+\frac{\gamma}{u-v}\right)\sum_{i=1}^{n-1}e_{ii}\otimes
e_{ii}+\left(1+\frac{\gamma}{u-v}\right)\frac{u(v+\gamma)}{v(u+\gamma)}\;e_{nn}\otimes
e_{nn}\nn\\
&&+\sum_{ij=1}^{n-1}e_{ii}\otimes
e_{jj}+\frac{(v+\gamma)}{v}\sum_{i=1}^{n-1}e_{ii}\otimes
e_{nn}+\frac{u}{(u+\gamma)}\sum_{i=1}^{n-1}e_{nn}\otimes
e_{ii}\nn\\
&&+\gamma\sum_{i=1}^{n}\sum_{k=1}^{i-1}\sum_{j=1}^{i-k}(e_{ii-k}\otimes
e_{jj+k-1}-e_{jj+k-1}\otimes
e_{ii-k})+\left(\frac{\gamma}{u-v}\right)\sum_{i\neq
j=1}^{n-1}e_{ij}\otimes
e_{ji}\nn\\
&&+\left(\frac{\gamma}{u-v}\right)\frac{u}{v}\sum_{i=1}^{n-1}e_{in}\otimes
e_{ni}+\left(\frac{\gamma}{u-v}\right)\frac{(v+\gamma)}{(u+\gamma)}\sum_{i=1}^{n-1}e_{ni}\otimes
e_{in}\nn\\
&&-\frac{\gamma}{(u+\gamma)}\sum_{i=1}^{n-1}\sum_{k=1}^{n-1-i}e_{nn-i}\otimes
e_{kk+i}+\frac{\gamma}{v}\sum_{i=1}^{n-1}\sum_{k=1}^{n-1-i}e_{kk+i}\otimes
e_{nn-i}.
\end{eqnarray}
It should be noticed that this $R$-matrix does depend on the two
spectral parameters $u$ and $v$ independently and not only through
their difference. This represents a novelty with respect to the
$R$-matrices which are solutions of the standard YBE, which usually
depend on the spectral parameters only through their difference.
Using \eqref{dynamicaltwistwithsp} and \eqref{dynamicalRmatrixsln},
a parametrization for the $D$-matrix in line with \eqref{paramD} is
realized. It represents an interesting example of decomposition for
non constant $D$-matrices.

\p Again, the matrix $Q$ \eqref{paramKbq} for the set of solutions
\eqref{Ksolutions} can be calculated. For instance, the $Q$ matrices
for the invertible solutions (\ref{Ksolutions}-\textbf{IIb}), which
do depend on the dynamical variables, are given by the following
triangular matrices
\begin{eqnarray}\label{Q}
\mbox{\textbf{IIb}}\qquad
Q=\sum_{i=1}^{(n-1)}\,e_{ii}+\left(1-\frac{n\,
f}{u}\right)\,e_{nn}+\sum_{i>j=1}^{(n-1)}\,\left(
  \begin{array}{c}
  i-1 \\
   i-j \\
   \end{array}
    \right)
f^{(i-j)}\,e_{ij}+\frac{M_{nj}}{u(u+\sigma-\gamma+f_0)} \,e_{nj},\nn\\
\end{eqnarray}
with
\begin{equation}
M_{nj}=\sum_{j=1}^{(n-1)}\,\left(
  \begin{array}{c}
  n-1 \\
   n-j \\
   \end{array}
    \right)
f^{(n-j)}u(u+\sigma)-\left(
  \begin{array}{c}
  n \\
   n-j+1 \\
   \end{array}
    \right)
f^{(n-j+1)}u+(-)^{(n-j+1)} n\,f
\sum_{\beta\in\mathcal{L}_{nj}}\,\prod_{i\in\beta}\;\lambda_i,\nn
\end{equation}
where $\mathcal{L}_{nj}$ is a set of set of indexes depending on $n$
and $j$: elements $\beta$ of $\mathcal{L}_{nj}$ are all distinct
sets of $(n-j+1)$ different indexes $k=1,\dots , n$. It can be
easily seen that \eqref{Q} collapses to \eqref{Qinfty} in the limit
$u\longrightarrow \infty$. This expression however does not allow
for a factorization of the dynamical shift in equation
\eqref{RQequation}, which, therefore, cannot be simplified to the
form \eqref{RQzeroequation}.

\section{New set of structure matrices for the SDRE}
\label{newsetofSM}

\p The work by Antonov et al. \cite{Antonov97} provides an
alternative parametrization for the matrix $^{trig}D$ corresponding
to the trigonometric case for which expression \eqref{DwithSP} is
the rational limit.\footnote{Note that recently, all these
parametrizations received a universal description in
\cite{Buffenoir07}.} This fact suggests the possibility to use such
a decomposition for finding new solutions for the matrices $A$, $B$
and $C$ satisfying the consistency equations
\eqref{ABCDconsitencyequations}. More precisely, the parametrization
provided in \cite{Antonov97} concerns the matrix $^{trig}D^T$ and
can be written as follows
\begin{equation}\label{newparametrizationforDT}
S_{12}\,\tilde{c}_1\,\tilde{c}_2(\lambda-\gamma
h_1)=\tilde{c}_2\,\tilde{c}_1(\lambda-\gamma
h_2)\,^{trig}D^{T}_{12};\quad S_{12}=d_1^{-1}\,\tilde{S}_{12}\,d_2,
\end{equation}
with
\begin{eqnarray}
&&^{trig}D_{12}\nn\\
&&\nn\\ &=&\sum_{i\neq
j=1}^n\left[\frac{\sinh(s+\gamma)}{\sinh(s)}\,e_{ii}\otimes
e_{ii}+\frac{\sinh(\gamma)\sinh(s+\lambda_{ij})}{\sinh(s)\sinh(\lambda_{ij})}\;e_{ij}
\otimes
e_{ji}+\frac{\sinh(\lambda_{ij}-\gamma)}{\sinh(\lambda_{ij})}\;e_{ii}
\otimes e_{jj}\right]\nn\\
\end{eqnarray}
\begin{eqnarray}
&&S_{12}\nn\\
&&\nn\\
&=&\sum_{i\neq
j=1}^n\left[\frac{\sinh(s+\gamma)}{\sinh(s)}\,e_{ii}\otimes e_{ii}
+\frac{\sinh(\gamma)}{\sinh(s)}\;e^{s(2(i-j)-n\,
\mbox{sign}(i-j))/n}\;e_{ij} \otimes e_{ji}+e^{\gamma\,
\mbox{sign}(i-j)}\;e_{ii} \otimes
e_{jj}\right]\nn\\
&&+2\sinh(\gamma)\left[\sum_{1=i<i'<j}^n\;e^{-2s(i'-i)/n}\;e_{ii'}
\otimes e_{jj'}- \sum_{i> i'>j=1}^ne^{-2s(i'-i)/n}\;e_{ii'} \otimes
e_{jj'}\right]\quad i+j=i'+j',\nn\\
\end{eqnarray}
where $s=(u-v)$ and the elements of the matrices $c$ and $d$ are
\begin{equation}
\tilde{c}_{jk}=e^{2 j(u+n \lambda_k)/n},\qquad d_{jk}=e^{2
j\gamma/n}\;\delta_{jk}.
\end{equation}
The notation adopted in writing these formulas has been adapted to
the present article, and therefore it differs slightly from the
conventions used in \cite{Antonov97}. The matrix $\tilde{S}$ from
\cite{Antonov97} is non-dynamical and depends on the spectral
parameters only through their difference $s$. According to Antonov
et al. it satisfies the YBE, and consequently, the matrix $S$
satisfies the YBE as well, since the following property holds
\begin{equation}
[S, d \otimes d]=0.
\end{equation}
As it is, the expression \eqref{newparametrizationforDT} implies a
parametrization for the $^{trig}D$-matrix which does not match the
decomposition \eqref{parametrizationforD}. However, it can be
noticed that both matrices $\tilde{S}$ and $^{trig}D$ are invariant
under the following transformation
\begin{equation}
\mbox{space 1}\longleftrightarrow \mbox{space 2}, \qquad
\gamma\longrightarrow -\gamma.
\end{equation}
This fact allows to rewrite expression
\eqref{newparametrizationforDT} as follows
\begin{equation}
S_{12}\,\tilde{c}_2,\tilde{c}_1(\lambda+\gamma
h_2)=\tilde{c}_1\,\tilde{c}_2(\lambda+\gamma
h_1)\,^{trig}D^{T}_{12},
\end{equation}
and consequently to obtain a parametrization for $^{trig}D$ in the
form \eqref{parametrizationforD}, namely
\begin{equation}\label{newparametrizationforD}
S^T_{12}\,c_1\,c_2(\lambda+\gamma h_1)=c_2\,c_1(\lambda+\gamma
h_2)\,^{trig}D_{12};\quad c_{j}=(\tilde{c}_j^{-1})^{T}\quad j=1,2.
\end{equation}
It can be noticed that unlike the case investigated previously in
section (\textbf{\ref{paramewithSP}}), the automorphism $g$ is set
equal to $\mathds{1}$, in spite of a spectral parameter dependence.
At this stage, putting in effect the procedure sketched in section
(\textbf{\ref{parametrization}}), it is possible to derive new
matrices $A$, $B$ and $C$ from the new set of data associated to
$^{trig}D$, namely $S^T=R$, $c=b$ and $g=\mathds{1}$. For instance,
one can think to obtain immediately the rational limit of the
matrices $c$ and $S^T$, hence suitable formulations for matrices $b$
and $R$ \eqref{paramD}, and consequently to find the corresponding
matrices $B$ and $A$ using expressions \eqref{paramB} and
\eqref{paramA}, respectively. Unfortunately, the rational limit of
matrix $c$, unlike the trigonometric case, leads to a matrix which
is non-invertible. An alternative possibility is first to find
matrices $B$ and $A$ in the trigonometric case, and subsequently to
take their well-defined, non-singular rational limit consistently
with equations \eqref{ABCDconsitencyequations}. This last procedure
turns out to be a better strategy. The new rational matrices $B$ and
$A$, which will be indicated as $\hat{B}$ and $\hat{A}$ to
differentiate them from the ACF matrices, are:
\begin{eqnarray}
\hat{B}_{12}&=&\sum_{i=1}^n\,e_{ii}\otimes p_{i}, \nn\\
p_i&=&\prod_{1=k\neq
i}^{n}\,\frac{\lambda_{ik}}{(\lambda_{ik}+\gamma)}\,e_{ii}+\sum_{1=j\neq
i}^n\left( e_{jj}-\gamma\prod_{1=k\neq
i,j}^{n}\,\frac{\lambda_{ik}}{\lambda_{jk}(\lambda_{ji}-\gamma)}\,e_{ij}
\right),\\
&&\nn\\
&&\nn\\
\hat{A}_{12}&=&\sum_{i=1}^n\,\left(1+\frac{\gamma}{s}\right)\,e_{ii}\otimes
e_{ii}+\sum_{i\neq
j=1}^n\,\left[\left(1-\frac{\gamma}{\lambda_{ij}}\right)\,e_{ii}\otimes
e_{jj}+\left(\frac{\gamma}{s}+\frac{\gamma}{\lambda_{ij}}\right)\,e_{ij}\otimes
e_{ji}\right]\nn\\
&&+\sum_{i\neq j=1}^n\,\prod_{k\neq i,j; l\neq
j}^{n}\frac{\gamma\lambda_{ik}}{\lambda_{jl}}\,\left(e_{ii}\otimes
e_{ij}-e_{ij}\otimes e_{ii}\right), \qquad s=(u-v).
\end{eqnarray}
Since matrix $B$ - called $\hat{B}$ in the present case - is known,
expression \eqref{paramB} can be used to find a suitable invertible
$b$ matrix, which, in the present case, turns out to be given by the
following expression
\begin{eqnarray}\label{newb}
\hat{b}_{ij}&=&\frac{\sum_{\alpha\in \mathcal{I}_{ij}}\,\prod_{l\in
\alpha}\,\lambda_{l}}{\prod_{k\neq j}^{n}\,\lambda_{jk}},
\end{eqnarray}
where $\mathcal{I}_{ij}$ is a set depending on $i$ and $j$. Each
element $\alpha$ of $\mathcal{I}_{ij}$ is a collection of $(n-i)$
different indexes $l\neq j$ and the total number of elements of this
set is given by the binomial coefficient $(n-1)!/(n-i)!(i-1)!$. For
instance, for $n=3$ \eqref{newb}becomes is:
\begin{equation}
\hat{b}=    \left(
       \begin{array}{cccc}
       \frac{\lambda_2\lambda_3}{\lambda_{12}\lambda_{13}} &
       \frac{\lambda_1\lambda_3}{\lambda_{21}\lambda_{23}} & \frac{\lambda_1\lambda_2}{\lambda_{31}\lambda_{32}}
       \\ \\
        \frac{\lambda_2+\lambda_3}{\lambda_{12}\lambda_{13}} & \frac{\lambda_1+\lambda_3}{\lambda_{21}\lambda_{23}} &
        \frac{\lambda_1+\lambda_2}{\lambda_{31}\lambda_{32}} \\ \\
        \frac{1}{\lambda_{12}\lambda_{13}}& \frac{1}{\lambda_{21}\lambda_{23}} & \frac{1}{\lambda_{31}\lambda_{32}}  \\
       \end{array}
     \right).
\end{equation}
Since matrices $\hat{b}$ and $\hat{A}$ are available, relation
\eqref{paramA} can be used for computing the $R$-matrix which
satisfies the ordinary YBE and which is
\begin{eqnarray}\label{newRmatrix}
\hat{R}_{12}&=&\sum_{i=1}^n\,\left(1+\frac{\gamma}{s}\right)\,e_{ii}\otimes
e_{ii}+\sum_{i\neq j=1}^n\,\left(e_{ii}\otimes
e_{jj}+\frac{\gamma}{s}\,e_{ij}\otimes
e_{ji}\right)\nn\\
&&+\gamma\sum_{i=1}^n\sum_{k=1}^{i-1}\sum_{j=1}^{i-k}(e_{i-ki}\otimes
e_{j+k-1j}-e_{j+k-1j}\otimes e_{i-k i}).
\end{eqnarray}
As expected, this matrix is non-dynamical. In addition it depends on
the spectral parameters only through their difference $s$. It can be
noticed that in the limit without spectral parameter, expression
\eqref{newRmatrix} becomes the transposed of matrix \eqref{Rinfty}.
Once again the rational Cremmer-Gervais-type matrix $\hat{R}$
exhibits a difference of one in root heights, compared to the
trigonometric case \cite{Buffenoir07}. This fact also explains why
one cannot take the direct trigonometric to rational limit in this
procedure, since the underlying non-dynamical matrices are
definitively of distinct form.

\p None can be said about the parametrization of the solutions $K$
of the SDRE with the structure matrices presented in this section,
since no $K$ matrices are known yet. However, it is interesting to
see whether the parametrization procedure could provide some
information concerning these unknown solutions and act as a shortcut
for finding them. Consider the situation with constant ACF structure
matrices, for which the `twist' is constructed from \eqref{binfty},
and consider also the parametrization for $D$
\eqref{parametrizationforD} with the `twist' built from
\eqref{newb}. Then, \eqref{RQequation} becomes
\begin{equation}\label{SDREQbq}
R_{12}\,Q_1\,\hat{b}_1\,Q_2(h_1)\,\hat{b}_1^{-1}
=Q_2\,\hat{b}_2\,Q_1(h_2)\,\hat{b}_2^{-1}\,\hat{R}_{12} \qquad
Q=b\,K\,\hat{b}^{-1}.
\end{equation}
All elements of this expression are known since the $K$ matrices
refer to solutions of the SDRE with ACF structure matrices. The same
cannot be said concerning the following expression
\begin{equation}\label{SDREQqb}
\hat{R}_{12}\,\hat{Q}_1\,b_1\,\hat{Q}_2(h_1)\,b_1^{-1}
=\hat{Q}_2\,b_2\,\hat{Q}_1(h_2)\,b_2^{-1}\,R_{12},\qquad
\hat{Q}=\hat{b}\,\hat{K}\,b^{-1},
\end{equation}
which is obtained from the SDRE using the new set of constant
structure matrices, for which the `twist' is constructed from
\eqref{newb} and the parametrization for the $D$-matrix is obtained
using the `twist' built from \eqref{binfty}. In this case the
matrices $\hat{K}$ are unknown. However, \eqref{SDREQbq} can be
manipulated in such a way to end up matching the formulation
\eqref{SDREQqb}. Writing \eqref{SDREQbq} as
\begin{equation}\label{SDREQbq-1}
\hat{R}_{12}\,[Q_1\,\hat{b}_1\,Q_2(h_1)\,\hat{b}_1^{-1}]^{-1}
=[Q_2\,\hat{b}_2\,Q_1(h_2)\,\hat{b}_2^{-1}]^{-1}\,R_{12},
\end{equation}
it can be verified that \eqref{SDREQqb} and \eqref{SDREQbq-1}
coincide provided
\begin{equation}
Q_2^{-1}(h_1)\,(\hat{b}_1^{-1}\,Q_1^{-1}\,b_1)=
(\hat{b}_1^{-1}\,\hat{Q}_1\,b_1)\,\hat{Q}_2(h_1),
\end{equation}
and therefore
\begin{equation}\label{relationsktildek}
(\hat{b}_2\,K_2^{-1}\,b_2^{-1})(h_1)\,K_1^{-1}=\hat{K}_1\,(\hat{b}_2\,\hat{K}_2\,b_2^{-1})(h_1),
\end{equation}
which represents a relation amongst the invertible constant
solutions $K$ and $\hat{K}$ of the SDRE with the two different sets
of structure matrices. Full investigation of this equation will be
left to future studies, however something can be said immediately
concerning the simplest case, namely $n=2$. In fact, making use of
the corresponding invertible set of solutions $K$
(\ref{KsolutionswoSP}-\textbf{IIb}), it is possible to compute
\begin{equation}\label{Qbqn2}
Q=\left(
                     \begin{array}{cc}
                       \phantom{-}0 & 1 \\
                       -1 & f+\sigma \\
                     \end{array}
                   \right).
\end{equation}
Though \eqref{Qbqn2} is dynamical, its dependence from the dynamical
variables appears through $\sigma$. This fact allows to simplify
\eqref{relationsktildek} which now reads
\begin{equation}\label{}
K_1^{-1}\,(\hat{b}_2\,K_2^{-1}\,b_2^{-1})(\sigma+\gamma)
=\hat{K}_1\,(\hat{b}_2\,\hat{K}_2\,b_2^{-1})(h_1).
\end{equation}
As a consequence
\begin{equation}
\hat{K}=K^{-1}=\left(
                   \begin{array}{cc}
                     1-\frac{f}{\lambda_{12}} & \frac{f}{\lambda_{12}} \\
                     -\frac{f}{\lambda_{12}} & \phantom{1}1+\frac{f}{\lambda_{12}} \\
                   \end{array}
                 \right),
\end{equation}
represents a set of invertible constant solutions of the SDRE for
the new set of constant matrices proposed in this section in the
case $n=2$. In the limit $f\longrightarrow 0$, the solution
$\hat{K}=\mathds{1}$ is obtained. More generally, when the mixed
matrices $Q$ \eqref{SDREQbq} built from the 2 cocycles $b$ and
$\hat{b}$ are quasi non-dynamical, $\hat{K}=K^{-1}$ always provides
a solution to the alternative SDRE. Unfortunately, in the present
case, for $n>2$ the dynamical dependence of matrices $Q$ cannot be
formulated in terms of $\sigma$ and therefore a more careful
investigation of equation \eqref{relationsktildek} is needed for
finding matrices $\hat{K}$.

\section{Conclusion}

\p The purpose of this article has been to extend previous work
concerning the classification of constant solutions of the SDRE to
the non-constant case, and to provide explicit realizations of the
parametrization procedure proposed in \cite{Avan07} for all elements
of the SDRE. During this analysis it has been shown how the
existence of two distinct parametrizations for the $D$-matrix leads
to different sets of structure matrices for the SDRE, and
consequently to new solutions $K$ for this equation. Because of the
parametrization procedure, it was possible to reveal a connection
amongst invertible solutions $K$ of the SDRE equation with the two
different sets of structure matrices available. In this context, an
explicit example has been provided for the case $n=2$. Still, a full
investigation of the SDRE equation is required for obtaining a
classification of the solutions $K$ related to the new set of
structure matrices. It will be interesting to see the relationship
amongst the integrable systems stemming from these solutions and the
RS models related to the ACF matrices, since all of them share the
same $D$-matrix. In addition, the existence, exemplified here, of
several inequivalent sets of $A$, $B$, $C$, $D$ matrices, which
share the same $D$-matrix with different de-twisting procedures
\eqref{parametrizationforD}, may explain why only two sets out of
four sets of constant solutions $K^{\infty}$ \eqref{KsolutionswoSP}
can be extended to the non-constant case \eqref{Ksolutions}. In
fact, besides the ACF set of structure matrices $A$, $B$, $C$, $D$
used in the present article, there may exist another set $A'$, $B'$,
$C'$, $D$ with the same limit $u,v,(u-v)\longrightarrow \infty$ and
different non-constant solutions $K'$, this time for the other two
sets of $K^{\infty}$. The new set $\hat{A}$, $\hat{B}$, $\hat{C}$,
$D$ found in section (\textbf{\ref{newsetofSM}}) does not realized
this scheme since $\hat{R}^{\infty}$, $\hat{b}^{\infty}$ are
different from $R^{\infty}$, $b^{\infty}$.

\p Key objects of the parametrization procedure are `twists' and
$R$-matrices for which explicit formulations are provided. Amongst
the $R$-matrices found, it is worth pointing out matrices
\eqref{R0matrixsln}, which satisfy a shifted YBE and which seem to
appear here for the first time. 
%Moreover, these matrices, together
%with \eqref{newRmatrix} and corresponding `twists' provide
%interesting examples of decomposition \eqref{parametrizationforD}
%for the $D$-matrix extended also to the non-constant case.
Existence of two decompositions of the matrix $D$ does not
contradict the uniqueness theorems in \cite{Buffenoir06,Buffenoir07}
since the decomposition \eqref{parametrizationforD} does not yield
an $R$-matrix solving the Yang-Baxter equation but its shifted
extension.

\p Building explicit monodromy matrices and consequently $N$-body
system or spin-chain Hamiltonians is now feasible. Once again, it
should be emphasized that in \cite{Avan07} it was shown how the
parametrization procedure proposed is able to extremely simplify
these constructions providing an elegant factorized form for the
monodromy matrices. This is due to the possibility to eliminate
completely the quantum-space shifts of the dynamical variables,
which are present in the original formulas \cite{Nagy05}, and which
make the construction of suitable monodromy matrices particularly
cumbersome.

\p Finally, it should be emphasized that the ACF matrices satisfying
the SDRE are associated to RS models in the bulk. In fact, the SDRE
is not a reflection equation in the usual sense and in
\cite{Arutyunov98} it was shown how any representation of the
algebra \eqref{SDRE} with the ACF structure matrices turns into a
representation of the fundamental relation $SLL=LLS$ provided
suitable transformations are applied. This fact can be also seen as
a rational and general consequence of the parametrization procedure
presented in section (\textbf{\ref{parametrization}}) and applied to
the ACF structure matrices. It would be interesting to study the RS
models with a boundary and to find a suitable algebra able to
describe them.

%\eject

\vskip 1.0cm \noindent{\bf Acknowledgements} \vskip .25cm \noindent
C.Z. thanks the \emph{Centre National de la Recherche Scientifique}
(CNRS) for the postdoctoral fellowship SPM 06-13, and J.A. thanks
\emph{Universidad do Algarve} for their hospitality.

\end{document}